\documentclass[J,12pt]{amsart}
\usepackage{epsfig,epsf,amsmath,amssymb}
\newtheorem{theorem}{Theorem}[section]

\newtheorem{definition}{Definition}
\newtheorem{example}[theorem]{Example}
\newtheorem{proposition}[theorem]{Proposition}
\newtheorem{remark}[theorem]{Remark}

\begin{document}
\title[]{Essential concepts of digital topology\\
(digital $k$-connectivity and $k$-adjacencis for digital products)}

\author[]{Sang-Eon Han}
\address[]{Department of Mathematics Education, Institute of Pure and Applied Mathematics\\
	Jeonbuk National University, Jeonju-City Jeonbuk, 54896, Republic of Korea\\
	e-mail address:sehan@jbnu.ac.kr, Tel: 82-63-270-4449.}

\begin{abstract}
		 The paper refers to several concepts which are essential to studying digital objects from the viewpoint of digital topology: digital $k$-connectivity or digital $k$-adjacency, $C$-compatible and normal $k$-adjacency for a digital product.
					Since L. Boxer has often mentioned the origins of these concepts  in an inaccurate way, we discuss something incorrectly cited or mentioned in Boxer's papers according to the facts. 
			\end{abstract}

%\date{2013. 3. 16}

\maketitle
\newpage

\section{\bf Terminology for studying digital images}\label{s1}

The present paper will often follow the notation.
${\mathbb N}$ indicates the set of natural numbers (or positive integers), ${\mathbb Z}$ means the set of  integers, $\lq\lq$$:=$" is used for introducing a new term and $\lq\lq$$\subset$" stands for a subset as usual. 
Besides,  for $a, b \in {\mathbb Z}$: the set of integers, 
$[a, b]_{\mathbb Z}:=\{t \in {\mathbb Z}\,\vert \,a \leq t \leq b\}$.
We also abbreviate the term $\lq\lq$digital topological" as $DT$.
In addition, $X^\sharp$  indicates the cardinality of the given set $X$.\\

\section{\bf Digital $k$-connectivities, digital $k$-neighborhood, a length between two points, and simple closed $k$-curves in ${\mathbb Z}^n, n \geq 2$}\label{s2}

In order to study digital images (or digital objects) from the viewpoint of digital topology, we had better consider digital $k$-connectivities (or digital $k$-adjacencies) for the $n$-dimensional lattice set ${\mathbb Z}^n, n \in {\mathbb N}$.
However, before 2004 (see \cite{H1,H2,H3}), focusing on the $3$-dimensional world, all papers in digital topology dealt with only  low dimensional digital image $(X, k)$, i.e., $X \subset {\mathbb Z}^3$ with one of the $k$-adjacency of
${\mathbb Z}^3$ \cite{B2,KR1,R1,R2}. 
For instance, the $2$-adjacency for ${\mathbb Z}$; the $4$- and $8$-adjacency for ${\mathbb Z}^2$; and the $6$-, $18$-, and $26$-adjacency for ${\mathbb Z}^3$ \cite{R1,R2}.\\
After generalizing these $k$-connectivities into those for $n$-dimensional cases, $n \in {\mathbb N}$, since 2004 (\cite{H1,H2,H3}), 
a digital image $(X, k)$ has often assumed to be a set $X$ in ${\mathbb Z}^n$ with one of the $k$-adjacency relations of
${\mathbb Z}^n$  from (2.2) below. 
Indeed, the $k$-adjacency relations for $X \subset {\mathbb Z}^n, n \in {\mathbb N}$, were initially established in \cite{H3} (see also \cite{H1,H2,H6,H11}).
 In detail, for a natural number $t$, $1 \leq t \leq n$, the distinct points
$p = (p_i)_{i \in [1,n]_{\mathbb Z}}$  and $q=(q_i)_{i \in [1,n]_{\mathbb Z}}\in {\mathbb Z}^n$
are $k(t, n)$-adjacent (see the page 74 of \cite{H3}) if 
$$\text{at most}\,\,t\,\,\text{of their coordinates  differ by}\,\,\pm1\,\,\text{and the others~coincide.}\eqno(2.1)$$

According to this criterion, the digital $k(t, n)$-connectivities of ${\mathbb Z}^n, n \in {\mathbb N}$, are formulated \cite{H3} (see also \cite{H11}) as follows:
$$k:=k(t,n)=\sum_{i=1}^{t} 2^{i}C_{i}^{n}, \text{where}\,\, C_i ^n:= {n!\over (n-i)!\ i!}. \eqno(2.2)$$
For instance, in ${\mathbb Z}^4$, 
$8$-, $32$-, $64$-, and $80$-adjacency; in ${\mathbb Z}^5$, 
$10$-, $50$-, $130$-, $210$-, and $242$-adjacency; and 
in ${\mathbb Z}^6$, 
$k(1,6)=12$, $k(2,6)=72$, $k(3,6)=232$, $k(4,6)=472$, $k(5,6)=664$, $k(6,6)=728$-adjacency
 are considered \cite{H1,H11}.\\
The first formula for the $k$-adjacency relations of
${\mathbb Z}^n, n \in {\mathbb N}$, is shown in \cite{H1,H2,H3} and the second version of them was suggested in \cite{H6} and the process of formulating them was mentioned in \cite{H8} (see Appendix of \cite{H8}). Eventually, the final version of them was established in \cite{H11} as in (2.2).\\
Besides, the notation $\lq\lq$$k(t, n)$" of (2.2) indeed gives us much information such as the dimension $n$ and the number $t$ associated with the statement of (2.1).
For instance, we may consider an $8$-adjacency in both ${\mathbb Z}^2$ and ${\mathbb Z}^4$.
Then, in the case that we follow the notation $k:=k(t,n)$, we can easily make a distinction between them by using  $k(2,2)=8=k(1,4)$ (see (2.1) and (2.2)).
Meanwhile, in Boxer and et al.'s  paper \cite{BK1}, the notation $c_2=8$ in  ${\mathbb Z}^2$ and $c_1=8$ in  ${\mathbb Z}^4$ are shown. 
Indeed, the notation of $c_u, u \in [1,n]_{\mathbb Z}$, is exactly based on the criterion of (2.1). However, Boxer did not refer to this origin. Boxer indeed took just only another notation with the same meaning as the criterion referred to in (2.1).
Even though the usage can be just a matter of preference, 
 Han's approach using $k:=k(t,n)$ can be helpful to study  more complicated cases such as digital products with a $C$-compatible, a normal, and a pseudo-normal $k$-adjacency \cite{H3,H7,H10,KHL1} and their applications.\\
Hereinafter, a relation set $(X, k)$ is assumed in ${\mathbb Z}^n$ with a $k$-adjacency of (2.2), called a digital image. Namely, we now abbreviate $\lq\lq$digital image $(X, k)$" as $\lq\lq$$(X, k)$" if there is no confusion.\\

Let us now recall a $\lq\lq$digital $k$-neighborhood of the point $x_0$" in $(X, k)$.
For $(X, k)$ and a point $x_0 \in X$, we have often used the notation
$$ N_k(x_0, 1):=\{x \in X \,\vert\,\,x\,\,\text{is}\,\,k\text{-adjacent to}\,\,x_0\} \cup \{x_0\}, \eqno(2.3)$$
which is called a digital $k$-neighborhood of the point $x_0$ in $(X, k)$.
The notion of (2.3) indeed facilitates some studies of digital images, e.g., digital covering spaces, digital homotopy, digital-topological ($DT$-, for brevity) $k$-group \cite{H15,H16}, and so on.
According to the literature on the digital $k$-neighborhood of (2.3), Kong and Rosenfeld \cite{KR1} and Bertrand \cite{B2} used the following notation (see the lines 14-16 from the bottom on the page 359 of \cite{KR1}), as follows:\\
In ${\mathbb Z}^2$, for $p \in {\mathbb Z}^2$,  
$N(p):=\{x\,\vert\,\,x\,\,\text{is}\,\,8\text{-adjacent to}\,\,p\} \cup \{p\}$ and\\
 in ${\mathbb Z}^3$, for $p \in {\mathbb Z}^3$,  $N(p):=\{x \,\vert\,\,x\,\,\text{is}\,\,26\text{-adjacent to}\,\,p\} \cup \{p\}$.\\
 In ${\mathbb Z}^3$, for $p \in {\mathbb Z}^3$,  $N^\ast(p):=\{x \,\vert\,\,x\,\,\text{is}\,\,26\text{-adjacent to}\,\,p\}$ \cite{BG1}.\\
To make a distinction from these notations, Han \cite{H1,H2,H3} used the notation of (2.3) because it has some advantages. 
Meanwhile, based on the notation of (2.3), we can represent it by using another term such as $N_k(x_0)$, $D_k(x_0)$, or $E_k(x_0)$ (see the page 3 of \cite{H12}), and so forth.
However, the notation of (2.3) has the following history. The notation of (2.3) indeed comes from the $\lq\lq$digital $k$-neighborhood of $x_0 \in X$ with
 radius $\varepsilon \in {\mathbb N}$" in \cite{H1,H2,H3} (see (2.4) below). 
 More precisely, for $(X,k), X \subset {\mathbb Z}^n$ and two distinct points $x$ and $x^\prime$ in $X$, a {\it simple $k$-path} from $x$ to $x^\prime$ with $l$ elements is said to be a sequence $(x_i)_{i \in [0, l]_{\mathbb Z}}$ in $X$ such that $x_0=x$ and $x_l=x^\prime$ and further, $x_i$ and $x_j$ are
$k$-adjacent if and only if $\vert i-j\vert=1$ \cite{KR1}. We say that a {\it length} of the simple $k$-path, denoted by
$l_k(x, y)$, is the number $l$ \cite{H3}.
For instance, let $(X=\{x_0=(0,0), x_1=(1,-1), x_2=(2,0), x_3=(1,1)\}, 8)$.
Then we observe that $l_8(x_0, x_1)=1, l_8(x_0, x_2)=2$, and $l_8(x_0, x_3)=1.$

\begin{definition} \cite{H1} (see also \cite{H3}) For a digital image $(X, k)$ on ${\mathbb Z}^n$,
	the digital $k$-neighborhood of $x_0 \in X$ with
	radius $\varepsilon$ is defined in $X$ to be the following subset of $X$
	$$N_k(x_0, \varepsilon) = \{x \in X \,\vert\, l_k(x_0, x) \leq
	\varepsilon\}\cup\{x_0\}, \eqno(2.4) $$
	where $l_k(x_0, x)$ is the length of a shortest simple $k$-path from $x_0$ to $x$ and  $\varepsilon\in {\mathbb N}$.
\end{definition}

For instance, let $Y=\{y_0=(0,0), y_1=(1,0), y_2=(1,1), y_3=(0,1), y_4=(-1,2)\}$ and 
$(Y,8)$. Then $N_8(y_0, 1)=Y \setminus \{y_4\}$	and consider $N_8(y_0, 2)=Y$.\\

Owing to Definition 1, it turns out that $N_k(x_0, 1)$ of (2.3) is a special case of $N_k(x_0, \varepsilon)$ of (2.4).
Note that the neighborhood $N_k(x_0, 2)$ plays an important role in establishing 
the notions of a radius $2$-local $(k_1,k_2)$-isomorphism and a radius $2$-$(k_1,k_2)$-covering map which is strongly associated with the
homotopy lifting theorem \cite{H2}  and is essential for calculating digital fundamental groups of digital images $(X,k)$ \cite{BK1,H3}. \\

We say that $(Y, k)$ is $k$-connected \cite{KR1} if for arbitrary distinct points $x, y \in Y$ there is a finite sequence $(y_i)_{i \in [0, l]_{\mathbb Z}}$ in $Y$ such that
$y_0=x$, $y_l=y$ and $y_i$ and $y_j$ are
$k$-adjacent if $\vert\, i-j \,\vert=1$ \cite{KR1}.
Besides, a singleton is assumed to be $k$-connected.\\

A simple closed $k$-curve (or $k$-cycle)
with $l$ elements in ${\mathbb Z}^n, n \geq 2$, denoted by $SC_k^{n,l}$ \cite{H3,KR1}, $l(\geq 4) \in {\mathbb N}$, is defined as the sequence $(y_i)_{i \in [0, l-1]_{\mathbb Z}}\subset {\mathbb Z}^n$ such that
$y_i$ and $y_j$ are $k$-adjacent if and only if  $\vert\, i-j \,\vert=\pm1(mod\,l)$.
Indeed, the number $l(\geq 4)$ is very important in studying $k$-connected digital images 
from the viewpoint of digital homotopy theory.
The number $l$ of $SC_k^{n,l}$ can be an even or an odd number. Namely, the number $l$ depends on the dimension $n$ and the $k$-adjacencies, as follows:

$$ \left \{
\aligned
& (1)\, \text{in the case of}\,\, k=2n (n \neq 2), \,\,\text{we have}\,\,l \in {\mathbb N}_0 \setminus \{2\};\\
& (2)\, \text{in the case of}\,\, k=4(n=2),\,\,\text{we obtain}\,\, l \in {\mathbb N}_0 \setminus \{2, 6\};\\
& (3)\, \text{in the case of}\,\, k=8(n=2),\,\,\text{we have}\,\,l \in {\mathbb N} \setminus \{1, 2, 3, 5\}; \\
& (4)\, \text{in the case of}\,\, k=18 (n=3),\,\,\text{we obtain}\,\,l \in  {\mathbb N} \setminus \{1, 2, 3, 5\};\,\text{and}\\
& (5)\, \text{in the case of}\,\, k:=k(t, n) \,\,\text{such that}\,\,3 \leq t \leq n,\\
&\,\text{we have}\,\,l \in {\mathbb N}\setminus \{1, 2, 3\}.
\endaligned\right\} \eqno(2.5) $$

\begin{remark}  (see (5) of \cite{H13}) The number $l$ of $SC_k^{n,l}$ depends on the dimension $n\in {\mathbb N}\setminus \{1\}$ and the number $l \in \{2a, 2a+1\}, a \in {\mathbb N} \setminus \{1\}$.
\end{remark}

\section{\bf $DT$-versions of the adjacencies of graph products established by Berge \cite{B1} and Harray \cite{H18}; a $C$-compatible and normal $k$-adjacency of a digital product established by Han \cite{H3,H7,H10}}\label{s3}

Based on the classical adjacencies of a Cartesian product of two graphs in \cite{B1} and \cite{H18}, the papers \cite{H3,H7,H10} introduced digital topological ($DT$-, for brevity) versions of them. More precisely, assume two digital images $(X, k_1)$ on ${\mathbb Z}^{n_1}$ and $(Y, k_2)$ on ${\mathbb Z}^{n_2}$ and further, consider a Cartesian product 
 $X \times Y \subset {\mathbb Z}^{n_1+n_2}$. Then we need some meaningful $k$-adjacencies of ${\mathbb Z}^{n_1+n_2}$ for $X \times Y$ to be a digital image $(X \times Y, k)$, where the $k$-adjacency of ${\mathbb Z}^{n_1+n_2}$  should be considered  according to (2.2).
Indeed, the adjacencies for a Cartesian product graph in \cite{H18} and a strong product graph in \cite{B1} facilitated the establishment of some $k$-adjacencis for a digital product $X \times Y \subset {\mathbb Z}^{n_1+n_2}$.
The former was helpful to develop the so-called $C$-compatible $k$-adjacency \cite{H10} (or the $L_C$-property of a Cartesian product in \cite{H7}) and the latter was motivated to develop a normal $k$-adjacency of $X \times Y$ \cite{H3}.
However, note that a $C$-compatible $k$-adjacency and a normal $k$-adjacency of $X \times Y$ need not be equal to the classical adjacencies of a Cartesian product graph in \cite{H18} and a strong product of graphs in \cite{B1}, respectively (see Example 3.5 and Remark 3.6).
More precisely, given two classical graphs, an adjacency of a Cartesian product graph in \cite{H18} always exists.
However, a $C$-compatible $k$-adjacency has its own features (see Example 3.5 and Remark 3.6).
Similarly, comparing with the adjacency of strong product graph in \cite{B1}, we clearly observe that a normal $k$-adjacency has its own features (see Example 3.5 and Remark 3.6).\\

 Motivated by the adjacency of a Cartesian product graph in \cite{H18}, a  $C$-compatible $k$-adjacency of a digital product was defined as follows:

\begin{definition} \cite{H3} ($DT$-version of the adjacency of a product graph in 
	\cite {H18})		
	For two digital images  $(X, k_1)$ on ${\mathbb Z}^{n_1}$ and $(Y, k_2)$ on ${\mathbb Z}^{n_2}$,
	consider a Cartesian product $X\times Y \subset {\mathbb Z}^{n_1+n_2}$. 
	We say that distinct two points $p:=(x, y)$ and $q:=(x^\prime, y^\prime)$ is $C$-compatible $k$-adjacent if and only if\\
	(1) $x=x^\prime$, $y$ is $k_2$-adjacent to $y^\prime$ or\\
	(2) $y=y^\prime$, $x$ is $k_1$-adjacent to $x^\prime$. \\
	Then we say that the $k$-adjacency of $(X\times Y, k)$ is $C$-compatible with the given two  $(X, k_1)$ and $(Y, k_2)$.	
\end{definition}

In Definition 2, even though the logic for establishing a $C$-compatible $k$-adjacency of $X \times Y$ seems to be similar to that for the adjacency of a Cartesian product graph in \cite{H18}, we need to note that the $C$-compatible $k$-adjacency of $X \times Y$ of Definition 2 is
different from the Cartesian
adjacency of \cite{H7} (see Example 3.5 and Remarks 3.6).

Owing to the digital $k$-neighborhood of (2.3), Definition 2 can be represented as follows:
\begin{definition} (A $C$-compatible $k$-adjacency of a digital product \cite{H10} and
	the $L_C$-property of a Cartesian product in \cite{H7}) \\
	Assume two digital images $(X_i, k_i)$ on ${\mathbb Z}^{n_i}, k_i:= k(t_i, n_i), i \in \{1, 2\}$. Then we say that a $k$-adjacency of $X_1\times X_2 \subset {\mathbb Z}^{n_1+n_2}$ is $C$-compatible with the given $k_i$-adjacency of $(X_i, k_i)$, $i \in \{1, 2\}$ if every point $p:=(x_1, x_2)$ in $X_1\times X_2$ satisfies the following property:
	$$N_k(p, 1)= (N_{k_1}(x_1, 1) \times \{x_2\})\cup (\{x_1\} \times N_{k_2}(x_2, 1)), \eqno(3.1)$$
	where the $k$-adjacency is one of the typical $k$-adjacencies of ${\mathbb Z}^{n_1+n_2}$ stated in (2.2).
\end{definition}

After comparing Definitions 2 and 3, the following is obtained.
\begin{proposition}\cite{H9} A $k$-adjacency of $X_1\times X_2 \subset {\mathbb Z}^{n_1+n_2}$ is $C$-compatible with the given $k_i$-adjacency of $(X_i, k_i)$, $i \in \{1, 2\}$, if and only if
every point $p:=(x_1, x_2)$ in $X_1\times X_2$ satisfies the following property:
$$N_k(p, 1)= (N_{k_1}(x_1, 1) \times \{x_2\})\cup (\{x_1\} \times N_{k_2}(x_2, 1)), \eqno(3.1)$$
where the $k$-adjacency is one of the typical $k$-adjacencies of ${\mathbb Z}^{n_1+n_2}$ stated in (2.2).
\end{proposition}

In Definitions 2 and 3, we strongly need to focus on the $k$-adjacency of $X\times Y \subset {\mathbb Z}^{n_1+n_2}$ as in (2.2). We can easily observe some difference between
a $C$-compatible $k$-adjacency and the adjacency of a graph product of \cite{H18}, as follows:

\begin{remark} Assume two digital images  $(X, k_1)$ on ${\mathbb Z}^{n_1}$ and $(Y, k_2)$ on ${\mathbb Z}^{n_2}$.
	Not every $X\times Y$ alway has a $C$-compatible $k$-adjacency in ${\mathbb Z}^{n_1+n_2}$ and further, the number $k$ of a $C$-compatible $k$-adjacency need not be unique (see Example 3.5 and Remark 3.6).
\end{remark}

Note that the notion of a $C$-compatible $k$-adjacency is essential to developing a 
$DT$-$k$-group (see \cite{H15,H16}).\\

Unlike the $C$-compatible $k$-adjacency of Definitions 2 and 3, motivated by the strong product graph in \cite{B1}, a  normal $k$-adjacency of digital product was defined as follows:

\begin{definition} \cite{H3} ($DT$-version of the adjacency of a strong product graph in \cite{B1})			
	For two digital images  $(X, k_1)$ on ${\mathbb Z}^{n_1}$ and $(Y, k_2)$ on ${\mathbb Z}^{n_2}$,
	consider a Cartesian product $X\times Y \subset {\mathbb Z}^{n_1+n_2}$. 
	We say that distinct two points $p:=(x, y)$ and $q:=(x^\prime, y^\prime)$ is normally 
	$k$-adjacent if and only if\\
	(1) $x=x^\prime$, $y$ is $k_2$-adjacent to $y^\prime$ or\\
	(2) $y=y^\prime$, $x$ is $k_1$-adjacent to $x^\prime$ or\\
	(3)	$x$ is $k_1$-adjacent to $x^\prime$ and $y$ is $k_2$-adjacent to $y^\prime$.\\
	Then we say that the $k$-adjacency of $(X\times Y, k)$ is normally $k$-adjacent with the given two  $(X, k_1)$ and $(Y, k_2)$
	\end{definition}

In Definition 4, we strongly focus on the $k$-adjacency of ${\mathbb Z}^{n_1+n_2}$ for $X\times Y$ according to (2.2). 
 In Definition 4, even though the logic for establishing a normal $k$-adjacency of $X \times Y$ seems to be similar to that for the adjacency of a strong product graph in \cite{B1}, it is clear that the normal $k$-adjacency of $X \times Y$ of Definition 4 is
different from the adjacency of \cite{B1} (see Example 3.5 and Remarks 3.6).

Using the digital $k$-neighborhood of (2.3), Definition 4 can be represented as follows:

\begin{definition} (A normal $k$-adjacency of a digital product \cite{H3,H10})
	Assume two digital images $(X_i, k_i)$ on ${\mathbb Z}^{n_i}, k_i:= k(t_i, n_i), i \in \{1, 2\}$. Then we say that a $k$-adjacency of $X_1\times X_2 \subset {\mathbb Z}^{n_1+n_2}$ is normal if
	every point $p=(x_1, x_2)$ in $X_1\times X_2$ satisfies the following property:
	$$N_k(p, 1)= N_{k_1}(x_1, 1) \times N_{k_2}(x_2, 1), \eqno(3.2)$$
	where the $k$-adjacency is one of the $k$-adjacencies of ${\mathbb Z}^{n_1+n_2}$ stated in (2.2).
\end{definition}
 Based on Definitions 4 and 5, the following is obtained.
 
 \begin{proposition} \cite{H10}
 	A $k$-adjacency of $X_1\times X_2 \subset {\mathbb Z}^{n_1+n_2}$ is normal  with the given $k_i$-adjacency of $(X_i, k_i)$, $i \in \{1, 2\}$ if and only if
 	every point $p=(x_1, x_2)$ in $X_1\times X_2$ satisfies the following property:
 	$$N_k(p, 1)= N_{k_1}(x_1, 1) \times N_{k_2}(x_2, 1), \eqno(3.2)$$
 	where the $k$-adjacency is one of the $k$-adjacencies of ${\mathbb Z}^{n_1+n_2}$ stated in (2.2).
 \end{proposition}

\begin{remark} Assume two digital images  $(X, k_1)$ on ${\mathbb Z}^{n_1}$ and $(Y, k_2)$ on ${\mathbb Z}^{n_2}$.
	Not every $X\times Y$ alway has a normal $k$-adjacency  in ${\mathbb Z}^{n_1+n_2}$ and further, a normal $k$-adjacency need not be unique (see Example 3.5 and Remark 3.6).
\end{remark}

 To confirm some difference between the adjacency of a Cartesian product graph in \cite{H18} and the $C$-compatible $k$-adjacency of a digital product and further, the adjacency of a strong product graph in \cite{B1} and the normal $k$-adjacency of a digital product of Definitions 4 and 5,
 we observe the following:
 
\begin{example} (1) $SC_4^{2,4} \times SC_4^{2,4}$ has a $C$-compatible $8$-adjacency.\\
	(2)  $SC_{26}^{3,5} \times SC_{26}^{3,5}$ has a $C$-compatible $k(3,6)$-adjacency and a normal $k(6,6)$-adjacency.\\
(3)	Each of $SC_8^{2,6} \times SC_8^{2,6}$  and  $SC_8^{2,4} \times SC_8^{2,6}$ has a normal $80$-adjacency.\\
(4)  $SC_8^{2,4} \times SC_8^{2,4}$ has a $C$-compatible $k$-adjacency, $k \in\{32=k(2,4), 64=k(3,4)\}$.\\
(5) Let $X=\{x_0=(0,0,0), x_1=(1,1,0), x_2=(1,2,1), x_3=(0,3,1), x_4=(-1,2,1), x_5=(-1,1,0)\}$ in ${\mathbb Z}^3$ with an $18$-adjacency. 
Indeed, $(X,18)$ is an example for $SC_{18}^{3,6}$.
Then $X^2$ has a normal $k$-adjacency, $k \in\{k(4,6)=472, k(5,6)=664, k(6,6)=728\}$.
\end{example}

\begin{remark} When comparing the adjacency of a Cartesian product graph in \cite{H18} and 
	a strong product graph in \cite{B1} and the $C$-compatible and the normal $k$-adjacency of a digital product, respectively, we observe some difference among them.
	Given two graphs, even though there are adjacencies for a Cartesian product graph in \cite{H18} and 
	a strong product graph in \cite{B1}, the $DT$-versions of them have some their own features, as follows:\\
	Let $MSC_{18}:=(x_i)_{i \in [0,5]_{\mathbb Z}}$ in ${\mathbb Z}^3$, where $x_0=(0,0,0), x_1=(1,-1,0), x_2=(1,-1,1), x_3=(2,0,1), x_4=(1,1 1), x_5=(1,1,0)$ as a special kind of an example for $SC_{18}^{3,6}$.
		Then the following is obtained.\\
		(1) $MSC_{18} \times MSC_{18} \subset {\mathbb Z}^6$ has neither a $C$-compatible $k$-adjacency nor a normal $k$-adjacency in ${\mathbb Z}^6$.\\
		(2) $SC_4^{2,4} \times SC_8^{2,6}$ does not have a normal $k$-adjacency in ${\mathbb Z}^4$, $k \in\{8,32,64,80\}$.
	\end{remark}

As mentioned in Example 3.5 and Remark 3.6, even though there are some big differences between the $C$-compatible $k$-adjacency of Han and the adjacency of Harray in \cite{H17}, and further, between the normal $k$-adjacency of Han and the adjacency of Berge in \cite{B1}, Boxer often does not make a distinction between them in his papers.

\begin{remark} (Utilities of the $C$-compatible and normal $k$-adjacency)\\
	(1) Calculation of digital fundamental groups of digital products \cite{H7,H10}.\\
	(2) The $C$-compatible $k$-adjacency is essential for developing $DT$-$k$-groups \cite{H15,H16}.
	\end{remark}

\section{\bf Further works}\label{s4}

We have discussed some concepts of digital topology which have been used to study digital topology.
Besides, we explained the $DT$-version of the $k$-adjacencies of digital products which incorrectly mentioned in Boxer's works.
Furthemore, as a further work, we will discuss some history of the notions of digital $k$-continuity, digital $k$-isomorphism, local $k$-isomorphism, and radius $2$-local $k$-isomorphism.
In addition,
based on the works \cite{B2,B3,BK1,H1,H2,H3,H5,H6,H14,H17,HP2,K1,P1,PZ1}, we will remark on Boxer's misconception on Han's papers \cite{H9,H14} and a revision of a pseudo-$(k_1, k_2)$-covering space.\\

{\bf Conflicts of Interest}: The author declares no conflict of interest.\\

\end{document}